\documentclass[a4paper,oneside,12pt]{article}
\usepackage{amssymb}

\newtheorem{theorem}{Theorem}[section]
\newtheorem{lemma}[theorem]{Lemma}

\newtheorem{proposition}[theorem]{Proposition}

\newtheorem{corollary}[theorem]{Corollary}

\newtheorem{hypothesis}[theorem]{Hypothesis \rm}

\newcommand{\ds}{\displaystyle}

\newcommand{\lra}{\longrightarrow}
\newcommand{\noi}{\noindent}

\newcommand{\St}{{\mathbf S}{\mathbf t}}
\def\matrice#1#2#3#4{\left( \begin{array}{cc} #1 & #2 \\
#3 & #4 \end{array}\right)}

\newcommand{\LL}{\mathbb L}
\newcommand{\PP}{\mathbb P}
\newcommand{\GG}{\mathbb G}

\newcommand{\CC}{\mathbb C}

\newcommand{\KK}{{\mathbb K}}
\newcommand{\bGG}{\bar {\mathbb G}}
\newcommand{\bLL}{\bar {\mathbb L}}
\newcommand{\bT}{\bar T}
\newcommand{\bt}{\bar t}
\newcommand{\HH}{\mathbb H}
\newcommand{\QQ}{\mathbb Q}

\newcommand{\Cinf}{{\mathcal C}_c^\infty}

\title{Branching laws for the Steinberg representation: the rank $1$ case}
\author{Paul Broussous}

\begin{document}
\maketitle
\tableofcontents

\vskip1.5cm

\noi {\bf Introduction}. Let $G/H$ be a reductive symmetric space over a $p$-adic field $F$, the algebraic
groups $G$ and $H$ being assumed semisimple of relative rank $1$. One of the  branching problems for the
Steinberg representation $\St_G$ of $G$ is the determination of the dimension of the intertwining space
${\rm Hom}_H (\St_G ,\pi )$, for any irreducible representation $\pi$ of $H$. In this work we do
not compute this dimension, but show how it is related to the dimensions of some other intertwining
spaces ${\rm Hom}_{K_i} ({\tilde \pi} ,1)$, for a certain finite family $K_i$, $i=1,...,r$, of anisotropic
subgroups of $H$ (here ${\tilde \pi}$ denote the contragredient representation, and $1$ the trivial
character).  In other words we show that there is a sort of ``reciprocity law'' relating two
different branching problems. 
\medskip

Before precisely stating our result, let us give a few motivations.  More generaly let $G$ be a locally
compact group and $H\subset G$ a closed subgroup. Fix an irreducible representation $\sigma$ of $G$
(belonging to a certain category of representations of $G$). The branching problem for $\sigma$
is the understanding of the restriction $\sigma_{\mid H}$ as a representation of $H$. In general this
restriction is not semisimple and a first task is to study the following two branching problems of
computing the dimensions of the intertwining spaces
$$
{\rm Hom}_H\, (\sigma ,\pi )\quad {\rm and} \quad {\rm Hom}_H\, (\pi ,\sigma )\ .
$$
\noi for all irreducible representations $\pi$ of $H$.  In the present work one is interested in the
first sort of intertwining space.  During the recent years there has been a tremendous progress in the
understanding of the intertwining spaces in the case where $G/H$ is a  symmetric space over
local field (cf. e.g. Beuzart-Plessis's talk \cite{BBki}
at Bourbaki seminar, and  Prasad's ICM 2018 contribution \cite{P}). 

Assume that $G/H$ is a Galois symmetric space over a non-archimedean local field $F$: there exist
a connected reductive group $\GG$ defined over $F$ and a Galois quadratic extension $K/F$
such that $G=\GG (K)$ and $H=\GG (F)$. In \cite{PReLa} Dipendra Prasad defines a quadratic character
$\omega_{K/F}$
of $H$ and states what he calls a {\it 'relative' Langlands conjecture} which, for an irreducible
representation $\sigma$ of $G$, computes
the dimension of ${\rm Hom}_H (\sigma,\omega_{K/F})$ in terms of a certain Galois parameter $\varphi_\sigma$
of $\sigma$. When $\sigma =\St_G$ is the Steinberg representation of $G$,  this conjecture coincides
with an earlier conjecture of Prasad (\cite{PSt} Conjecture 3). When $F$ has characteristic $0$, this
latter conjecture was proved  in
\cite{BP} by Beuzart-Plessis. Particular cases of this conjecture where also proved by the author and by
F. Court\`es (\cite{Br}, \cite{Cou}, \cite{Cou2}), when the characteristic of $F$ is odd,  and
with some more explicit features
(explicit intertwining operators and explicit test vectors).
\medskip

In the present work we tackle a more general branching problem but in a very modest framework: $G/H$ is a
reductive symmetric pair (not necessarily Galois), where $G$ and $H$ are semisimple of relative rank
$1$ over $F$, $\sigma =\St_G$, the Steinberg representation of $G$ and  $\pi$ is any irreducible
representation of $H$.  So fix such a reductive symmetric pair $G/H$ and assume that it is
attached to some rational involution $\theta$ of $G$. It is a theorem of Helminck and Wang that $H$
has a finite number of orbits in the flag variety of $G$, that is the set of minimal parabolic $F$-groups
of $G$. We denote by $r$ the number of orbits of $H$ in the set of minimal parabolic $F$-groups
that are not $\theta$-stable. Fix a set of representatives $P_1$, ..., $P_r$ of those orbits, and set
$K_i =H\cap P_i$. We prove that the $K_i$'s are anisotropic subgroups of $H$. Our main result is:
\medskip

\noi {\bf Theorem 1}. {\it  Let $\St_H$ denote the Steinberg representation of $H$. Let $\pi$ be
an irreducible smooth representation of $G$ assume to be non isomorphic to the trivial character of $H$.
Then we have 
  $$
   {\rm dim}_\CC \, (\St_G ,\pi )=
          \left\{\begin{array}{cc}
              {\ds  \sum_{i=1,...,r}} {\rm dim}_\CC \, {\rm Hom}_{K_i}({\tilde \pi}, 1) &
              {\rm if}\ \ \pi\not\simeq \St_H\\ 
              {\ds  \sum_{i=1,...,r}} {\rm dim}_\CC \,   {\rm Hom}_{K_i}({\tilde \pi}, 1)+1 & {\rm if}\ \ 
              \pi\simeq\St_H
          \end{array}\right.
          $$}

As a byproduct of the proof of Theorem $1$, we obtain certains  values of the derived functor ${\rm Ext}^*_H$.
\medskip

\noi {\bf Theorem 2}. a) {\it Let $\pi$ be an irreducible smooth representation of $\pi$, not isomorphic
  to the trivial character of $H$. Then ${\rm Ext}_H^k (\St_G ,\pi )=0$, for all $k\geqslant 1$. }

\noi b) {\it We have ${\rm Ext}_H^k (\St_G ,1_H )=0$, for all $k\geqslant 2$. Moreover
  ${\rm EP}_H (\St_G ,1_H )=r-1$, where ${\rm EP}$ denotes the Euler-Poincar\'e characteristic of
  ${\rm Ext}_H^*$.}
\medskip

Somehow this work raises more questions than gives answers. For instance, it would be nice to have an
estimate of the number $r$ (in all our  calculations we find $r\leqslant 2$). Similarly it would be good
to know whether there are example giving ${\rm Ext}_H^1 (\St_G ,1_H )\not= 0$ or not. Finally the quotients
$H/K_i$, $i=1,...,r$,  are often proved to be reductive symmetric spaces. Is that always true?  Note that
rank $1$ reductive groups are explicitely described: their all are isogeneous to (forms of) classical
groups (e.g. see \cite{Ca}). So one could imagine to handle these questions by a case by case procedure.
However to the author's knowledge there is no known classification
of rational involutions for that class of groups. 
\medskip

Since all objects and assertions of this work make sense when $F$ is a finite field, we decided to include
that case and give an uniform treatment.
\medskip

\noi {\bf Outline}. The notation is introduced in {\S}1. In {\S2} we specialize the general results of
Helminck and Wang on rationality properties of involutions \cite{HW} to the rank $1$ case. In particular
we introduce the anisotropic groups $K_i$. Our main Theorem 1 and 2 are proved in {\S}3. In {\S}4 we
entirely work out the particular case of the Galois  pair associated to ${\rm PGL}(2)$.  Open questions
are discussed in {\S}5.

\medskip

\noi {\bf Acknowledgements}. The reading of \cite{P} was a great source of inspiration for writing this paper.
I must thank Dipendra Prasad, Abderrazak Bouaziz and Nadir Matringe for numerous conversations during
the elaboration of this work (in particular I thank Prasad and Bouaziz for their help in  the proof of
Proposition \ref{symspace}).

\section{Notation and hypotheses}
Throughout this article we use the following notation:
\smallskip

$F$ is a either a non-archimedean, non-discrete, locally compact
field\footnote{Such a field will be called ``{\it $p$-adic}''}, or a
  finite field. In any case we assume $F$ to have {\it odd characteristic}.  

  If $\mathbb V$ is an algebraic variety defined over a field $k$, we denote
  by ${\mathbb V}(k)$, or by ${\mathbb V}_k$, its set of $k$-rational points. 
  
$\mathbb G$ is a connected semisimple algebraic group defined over
$F$, assumed to have relative rank $1$,

$\theta$ is an $F$-rational involution of $\mathbb G$,

$\mathbb H$ denotes  the fixed point group
${\mathbb G}^\theta$ and ${\mathbb H}^\circ$ its connected component. 
\smallskip

 The $F$-algebraic group $\mathbb H$ is always reductive. We make the
 following assumption:

 \begin{hypothesis} The group $\mathbb H$ is semisimple of relative
   rank $1$.
 \end{hypothesis}

 Let us give some examples.
 \medskip

 \noi {\it Example 1}. Here $F$ is a $p$-adic field.  Let $D$ be a central quaternion $F$-algebra.
 Fix an element $\tau\in D$ such that $\tau^2 \in F\backslash
 F^2$ and write $E$ for the field $F[\tau ]\subset D$. Let $\mathbb G$ be
 the semisimple $F$-group attached to ${\rm SL}(2,D)$. The inner automorphism
 ${\rm Int}(\tau I_2 )$, where $I_2$ denotes the identity matrix,
 is an $F$-rational involution of
 $\mathbb G$. The group ${\mathbb H}={\mathbb G}^\theta$ is the connected
 semisimple $F$-group attached to ${\rm SL}(2,E)$, that is the restriction of scalars
 ${\rm Res}_{E/F} \, {\mathbb S}{\mathbb L}_2$.
 
 \smallskip

 \noi {\it Example 2}. Let $F$, $D$, $\tau$ and $E$ be as above and let
 $\mathbb G$ be the semisimple $F$-group attached to ${\rm
   PGL}(2,D)$. Let $\theta$ be the rational involution of $\mathbb G$
 given by the inner automorphism ${\rm Int}( \tau I_2 )$. Let
 $\sigma$ be the generator of ${\rm Gal}(E/F)$. Then by
 Skoelem-Noether's theorem, the natural action of $\sigma$ on ${\rm
   M}(2,E)\subset {\rm M}(2,D)$ is induced by ${\rm Int}(x_\sigma )$ for
 some element $x_\sigma\in {\rm GL}(2,D)$. Then $\mathbb H$ is the
 non-connected semisimple algebraic $F$-group attached to ${\rm
   PGL}(2,E)\ltimes \{ 1, {\bar x}_\sigma\}$, where ${\bar x}_\sigma$
   is the image of $x_\sigma$ in ${\rm PGL}(2,D)$.  
   \medskip

\noi {\it Example 3}. {\it The Galois Case}. Let $\mathbb H$ be a connected semisimple
$F$-group of relative rank $1$ and  ${\rm Res}_{E/F}{\mathbb
  H}/{\mathbb H}$ be the associate {\it Galois symmetric space},  for some
Galois quadratic extension $E/F$. Then if  the $E$-rank of $\mathbb H$
is $1$, the groups ${\mathbb G}= {\rm Res}_{E/F}{\mathbb H}$ and
${\mathbb H}$ satisfy Hypothesis 1.1. 
\medskip

\noi {\it Example 4}. Here $F$ is a $p$-adic field.
Let $(V,q)$ be a quadratic space of rank $5$ over $F$ such
that $V$ writes $V_0 \perp H$, where $(V_0 ,q_{\mid V_0})$ is
anisotropic and $(H,q_{\mid H})$ is a hyperbolic plane. Let $v_1$ be
an anisotropic vector of $H$ and set  $W=v_1^\perp$. Assume than $W$
is an  isotropic subspace of $V$. Then $G={\rm SO}(V,q)$ and $H={\rm
  SO}(W,q_{\rm \mid W})$  are
semisimple groups of $F$-rank $1$. Let $s\in {\rm O}(V,q)$ be the
reflection of hyperplane $W$. The the action of $s$ on $G$ by
conjugation induces an $F$-rational involution $\theta$  and
we have $H=(G^\theta )^\circ$.  

\medskip

 We set $H={\mathbb H}^\circ (F)$ and $G={\mathbb G}(F)$. These groups are
 locally compact and totally disconnected. We denote by $\St_G$
 the Steinberg representation of $G$. It is defined as
 follows. Let $\mathbb P$ be a parabolic subgroup of $\mathbb G$ defined
 over $F$ and minimal for this property; set $P={\mathbb P}(F)$. Let
 $V$ (resp. $V_0$) be the space of locally constant (resp. constant)
 complex valued functions on $G/P$. Then $\St_G$ is the
 irreducible representation $V/V_0$.

\section{Orbits in the flag variety} In \cite{CLT}, Curtis, Lehrer and
Tits attach to any connected reductive group $\Gamma$ defined
over a field $k$ a topological space ${\mathcal B}(\Gamma,k)$ called the
(Curtis-Lehrer-Tits) spherical building. It enjoys the following
functorial property: attached to any monomorphism of connected
$k$-reductive groups $f$~: $\Gamma_1 \longrightarrow \Gamma_2$,
there is a continuous embedding ${\mathcal B}(f)$~: ${\mathcal
  B}(\Gamma_1 ,k)\lra {\mathcal B}(\Gamma_2 , k)$. Moreover this embedding is
$\Gamma_1 (k)$-equivariant.
\smallskip

 When $\Gamma$ is semisimple of $k$-rank $1$, then as a $\Gamma
 (k)$-set, ${\mathbb B}(\Gamma , k)$ is isomorphic to the $\Gamma
 (k)$-set of minimal  parabolic $k$-subgroups of $\Gamma$. 
 \smallskip

 A remarkable result is the following.

 \begin{theorem} (\cite{KM} Theorem 3.3.1) Assume that ${\rm
     char}(k)\not= 2$, and let $\theta\in {\rm Aut}_k (\Gamma )$ be a
   rational involution. Set $\Gamma '=(\Gamma^\tau )^\circ$ and let
   $i$~: $\Gamma ' \lra \Gamma$ be the natural embedding. Then
   the image of ${\mathcal B}(i)$~: ${\mathcal B}(\Gamma ',k)\lra
   {\mathcal B}(\Gamma ,k)$ is  the fixed point set  ${\mathcal B}(\Gamma ,k)^\tau$.
 \end{theorem}

 Let $\mathbb G$ be as in the introduction. We let ${\mathcal
   P}_F ({\mathbb G})$ (resp. ${\mathcal P}_F ({\mathbb H})$) denote
 the set of minimal (that is proper)   parabolic $F$-subgroups of $\mathbb G$
 (resp. of $\mathbb H$). As a corollary of the previous theorem, we
 have the following.

 \begin{corollary} a) We have a $H$-equivariant injection
   $$
   {\mathcal P}_F ({\mathbb H})\hookrightarrow {\mathcal
     P}_F ({\mathbb G})
   $$
\noi whose image consists of those  minimal parabolic $F$-subgroups of
$\mathbb G$ which are $\theta$-stable.

\noi b) With the notation of a), if $\mathbb Q$ is a minimal parabolic
subgroup of $\mathbb H$ which is mapped to ${\mathbb P}\in {\mathcal
     P}_F ({\mathbb G})^\theta$, then ${\mathbb P}\cap \HH^\circ =
\QQ$ and $\PP (F)\cap H =\QQ (F)$.
 \end{corollary}

We now collect some important  facts on parabolic subgroups and rational involutions. The
material is taken from \cite{HW}.
\smallskip

 A parabolic $F$-subgroup $\PP$ of $\GG$ is called $\theta$-split if
 $\PP^\theta := \theta (\PP )$ and $\PP$ are opposite parabolic
 subgroups. Since $\GG$ is of $F$-rank $1$, a parabolic
 $F$-subgroup of $\GG$ is either $\theta$-split or $\theta$-stable, indeed in that
 case two different parabolic $F$-subgroups are always opposite.
 
 \smallskip

  Let $A$ be a $\theta$-stable torus of $\GG$. We define $A^+$ and
  $A^-$ to be the maximal subtori of $A$ such that $\theta_{\mid A^+}$
  is the identity and $\theta_{\mid A^-}$ acts as $t\mapsto t^{-1}$. A
  $\theta$-stable $F$-torus $S$ of $G$ is called $(\theta ,F)$-split if it is
  $F$-split and if $S=S^-$.  
  \smallskip

  We shall need the following result.

  \begin{proposition}\label{HW} (\cite{HW} Lemma 2.4 and Proposition 4.7.) a)
    Any minimal parabolic $F$-subgroup of $\GG$ contains a
    $\theta$-stable maximal $F$-split torus of $\GG$.

    \noi b) Let $P$ be a
    $\theta$-split parabolic $F$-subgroup of $\GG$ and $T$ a
    $\theta$-stable maximal $F$-split torus of $P$. The following
    assertions are equivalent:
    \smallskip

    (i) $P$ is a minimal $\theta$-split parabolic $F$-subgroup of
    $\GG$.

    (ii) $T$ is a maximal $(\theta ,F)$-split torus of $G$ and $P\cap
    \theta (P) = {\rm Z}_G (T^- )$.
  \end{proposition}

   Since our group is of $F$-rank $1$ all non-trivial $\theta$-split
   and $F$-split $F$-tori are maximal $(\theta ,F)$-split tori and all
  $\theta$-split proper parabolic $F$-subgroups of $\GG$ are minimal
   $\theta$-split parabolic $F$-subgroups of $\GG$. Therefore we have
   the following: 

   \begin{corollary} Let $P$ be a $\theta$-split proper parabolic
     $F$-subgroup of $\GG$. Then $P$ contains a maximal $F$-split torus
     $T$ of $\GG$ 
     which is $\theta$-split and a Levi decomposition of $P$ is
     $P=MU$, where $U$ is the unipotent radical of $P$ and
     $M=Z_{G}(T)$.
   \end{corollary}

  \begin{proposition} Let $P$ be a $\theta$-split minimal parabolic
    subgroup of $\GG$. Then $P\cap \HH$ is an $F$-anisotropic subgroup of
    $\GG$. In particular if $F$ is a local field, $(P\cap \HH)_F$ is a
    compact subgroup of $H = \HH^\circ_F$.
  \end{proposition}

  \noi {\it Proof}. Let $T\subset P$ be a rank $1$ $(\theta ,
  F)$-split torus, $M=Z_G (T)=P\cap \theta (P )$. The Levi subgroup
  $M$ writes as an almost direct product $TM_a$: $T$ is the maximal
  central $F$-split torus of $M$ and $M_a$ is $F$-anisotropic. In
  particular $T\cap M_a$ is finite.

   We have $P\cap \HH = (P\cap \theta (P))\cap \HH = M\cap \HH$. Let
   $S\subset M\cap\HH$ be an $F$-split torus. Then $S$ is necessarily
   central in $M$, otherwise $M$ would have positive semisimple
   $F$-rank. It follows that $S\subset T$. So $S\subset T\cap H$ $=$
   $\{ t\in T \ ; \ t^2 =1\}$ $=$ $\{\pm 1\}$ and $S$ is
   trivial. Hence $P\cap M$ is anisotropic. 
   \medskip

   We recall the following striking fact.

   \begin{theorem} \label{HWfinite} (\cite{HW} Corollary 6.16.) The set of
     $H$-orbits in
     ${\mathcal P}_F (\GG )$ is finite.
   \end{theorem}

 By Corollary 2.2, one may identify  the subset ${\mathcal P}_F (\HH )\subset {\mathcal
   P}_F (\GG )$ with ${\mathcal P}_F (\GG )^\theta$. Hence the  $\theta$-stable minimal parabolic
 $F$-subgroups of $\GG$ form a single $H$-orbit. We denote by $r$ the number of
 $H$-orbits in the set ${\mathcal P}_F (\GG )\backslash   {\mathcal
   P}_F (\HH )$ of minimal $\theta$-split parabolic $F$-subgroups of $G$.
 \medskip

 \noi {\bf Examples}. Let $E/F$ be a quadratic extension of
 non-archimedean local fields. When $\GG /\HH = {\rm Res}_{E/F} {\rm
   PGL}(2)/{\rm PGL}(2)$, we have $r=1$. When $\GG /\HH ={\rm
   Res}_{E/F}{\rm SL}(2) /{\rm SL}(2)$, we have $r=2$.  In Example 2
 of {\S}1, we have $r=1$. In Example 4 of {\S}1, we have $r\leqslant
 2$. Proofs are left to the reader as exercises. 
 \medskip

 We fix once for all a set $\{ P_1 ,P_2 ,...,P_r  \}$ of
 representatives of the $H$-orbits of $\theta$-split parabolic
 $F$-subgroups of $\GG$. For $i=1,...,r$, let us  denote by
 $\Omega_i$ the orbit $H .P_i$ and set  ${\mathbb K}_i = P_i \cap
 \HH$, $K_i ={\mathbb K}_i (F)$.

 \begin{proposition} Assume that $F$ is a local non archimedean
   field so that ${\mathcal P}_F (\HH )$ and ${\mathcal P}_F (\GG )$
   are naturally compact, totally disconnected, topological spaces. 
   \smallskip

   \noi a) The subset of  ${\mathcal P}_H (\GG )$ formed of $\theta$-stable
   parabolic subgroups is homemorphic to ${\mathcal P}_F (\HH )$.

   \noi b) For $i=1,...,r$ equip $\Omega_i$ with the induced topology
   and $\HH_F /K_i$ with the quotient topology.
   \smallskip

   (i) For $i=1,...,r$, $\Omega_i$ is open in ${\mathcal P}_F (\GG )$. 
   
   (ii) For $i=1,...,r$, $\Omega_i$ and
   $H /K_i$ are homeomorphic.
 \end{proposition}

 {\it Proof}.  Helminck and Wang call a minimal parabolic $F$-subgroup of
 $G$ {\it quasi $\theta$-split} if $P$ is contained in a minimal
 $\theta$-split parabolic $F$-subgroup of $G$.  So in our case,
 for $i=1,...,r$, $P_i$ is quasi $\theta$-split. By  \cite{HW} Proposition
 13.4, if $P$ is a minimal parabolic $F$-subgroup of $G$, $\HH_F .P_F$ is open
 in $G$ if, and only if, $P$ is quasi $\theta$-split. Hence point b) (i)
 follows. The other assertions of the Proposition are standard.

 \section{Branching laws and their Ext versions}

 Let ${\mathcal S}(H)$ be the abelian category of smooth complex representations of $H$. We denote by
 ${\rm Ext}_H^k$, $k\geqslant 0$ the derived functors of ${\rm Hom}_H$.
 
 In this section we investigate the $\CC$-vector  spaces
 ${\rm Ext}^k_H \, (\St_G, \pi )$, $k\geqslant 0$, for any irreducible smooth
 complex representation $\pi$ of $H$. 
 When the $\CC$-vector spaces ${\rm Ext}_H^k \, (\St_G ,\pi )$ are finite dimensional and vanish for $k$ large,
 the Euler-Poincar\'e caracteristics
 $$
 {\rm EP}_H \, (\St_G ,\pi ):=\sum_{k\geqslant 0} (-1)^k {\rm dim}_\CC \, {\rm Ext}^k\,
 (\St_G ,\pi )
 $$
 \noi is well defined.
 \bigskip

 In the sequel we abreviate $X={\mathcal P}_F (\GG )$ and
 $Y={\mathcal P}_F(\HH )$. We identify $Y$ with a closed subset of $X$ in a
 $H$-equivariant way. The open subset $\Omega := X\backslash Y$ writes $\Omega
 =\Omega_1 \sqcup \cdots \sqcup \Omega_r$ as in {\S}2. 

 For any topological space $Z$, write $\Cinf (Z)$ for the ${\mathbb C}$-vector
 space of locally constant functions $f$~: $Z\lra {\mathbb C}$ with compact
 support. If $T\subset Z$ if an open subset, we have a natural linear map
 $i_{T,Z}$~: $\Cinf (T)\lra \Cinf (Z)$ extending functions by zero. Similarly if
 $T\subset Z$ is a closed subset, we have the natural restriction map
 $r_{Z,T}$~: $\Cinf (Z)\lra \Cinf (T)$.

 We need the following well known fact due to Bernstein and Zelevinsky
(\cite{BZ}, Proposition 1.8). 

 \begin{lemma} Let $Z$ be a totally disconnected locally compact topological
   space and $U$ be an open subset of $Z$. Write $T=Z\backslash U$. Then the
   maps $i_{U,Z}$ and $p_{Z,T}$ fit into a short exact sequence:
   $$
   0\lra \Cinf (U)\lra \Cinf (Z)\lra \Cinf (T)\lra 0\ .
   $$stable.
 \end{lemma}

 If $\Gamma_1$ is a locally compact totally disconnected topological group, 
 $\Gamma_2$ is a closed subgroup of $\Gamma_1$, and $\sigma$ is a smooth
 complex representation of $\Gamma_2$, we denote by
 ${\rm ind}_{\Gamma_1}^{\Gamma_2}\, \sigma$ the compactly smooth induced
 representation (see \cite{BZ}, Definitions 2.21, 2.22 for the definition of the basic induction
 functors).

   Fix a $\theta$-stable minimal parabolic $F$-subgroup $P$ of $\GG$, so that
   $Q:= P\cap \HH$ is a  parabolic $F$-subgroup of $\HH$. We
   observe that the $G$ representations $\Cinf (X)$ and ${\rm ind}_{P_F}^G \, 1$
   are isomorphic (here $1$ denotes the trivial representation of $P_F$), and
   that the $H$-representations $\Cinf (Y)$ and ${\rm ind}_{Q_F}^H \, 1$ are
   isomorphic as well. Similarly, it follows from Proposition 2.7.b(ii)
   that for $i=1,...,r$,
   the $H$-representations $\Cinf (\Omega_i )$ and ${\rm ind}_{K_i}^H \, 1$
   are isomorphic.

   \begin{proposition} a) We have an isomorphism of $H$-modules
     $$
     \Cinf (\Omega ) \simeq \bigoplus_{i=1,...,r} {\rm ind}_{K_i}^H \, 1 \ .
     $$
     \noi b) We have a short exact sequence of $H$-modules:
     $$
     0\lra \bigoplus_{i=1,...,r} {\rm ind}_{K_i}^H \, 1 \ \lra \St_G \lra
     \St_H \lra 0\ .
     $$
     \end{proposition}

   \noi {\it Proof}. Applying Lemma 3.1 to $Z=X$ and $U=\Omega$, we get the
   short exact sequence of $H$-modules:
   $$
   0\lra \Cinf (\Omega )\lra {\rm ind}_{P_F}^G \, 1\lra {\rm ind}_{Q_F}^H \, 1
   \lra 0
   $$

   \noi Moreover if $\CC_X$ and $\CC_Y$ denote the space of constant complex
   functions on $X$ and $Y$ respectively, and $0_\Omega$ denote the space of zero
   functions on $\Omega$, the previous exact sequence restricts to
   $$
   0\lra 0_\Omega \lra \CC_X \lra \CC_Y \lra 0
   $$
   \noi By dividing, we obtain the short exact sequence of $H$-modules:
   $$
   0\lra \Cinf (\Omega )\lra \St_G\lra \St_H \lra 0\ .
   $$

   \noi We are thus reduced to proving a). For this we prove by induction on
   $u=1,...,r$, that
   $\Cinf (\Omega_1 \sqcup \cdots \sqcup \Omega_u )={\rm ind}_{K_1}^H
   \, 1  \oplus \cdots \oplus {\rm ind}_{K_u}^H \, 1$. This is true for $u=1$.
   Assume this is true for $u<r$. Apply Lemma 3.1 to
   $X=\Omega_1 \sqcup \cdots    \sqcup \Omega_{u+1}$ and
   $U=\Omega_1 \sqcup \cdots \sqcup \Omega_u$. We obtain
   the short exact sequence of $H$-modules:
   $$
   0\lra \Cinf (\Omega_1 \sqcup \cdots \sqcup \Omega_u )\lra
   \Cinf (\Omega_1 \sqcup \cdots \sqcup \Omega_{u+1}) \lra {\rm ind}_{K_{u+1}}^H \,
   1\lra 0
   $$
   \noi Since $H$ is semisimple and $K_{u+1}$ compact, the representation
   $ {\rm ind}_{K_{u+1}}^H \, 1$ is a projective object in the category of smooth
   complex representations of $H$. Hence the previous short exact sequence
   splits and we are done.
   \bigskip

   From now on we fix an irreducible  smooth complex representation $\pi$ of
   $H$. Applying the functor ${\rm Hom}_H \,  (-, \pi )$,  from the category of
   smooth complex representations of $H$ to the category of $\CC$-vector spaces,
   to the short exact sequence b) of
   Proposition 3.2, we obtain the cohomology long exact sequence:
   $$
   \begin{array}{ccccccc}
     0 & \lra & {\rm Hom}_H\, (\St_H ,\pi ) & \lra & {\rm Hom}_H (\St_G ,\pi ) &
     \lra & \ds \bigoplus_{i=1,...,r} {\rm Hom}_H \,  ({\rm ind}_{K_i}^H \, 1 , \pi )\\
     &  \lra & {\rm Ext}^1_H\, (\St_H ,\pi ) & \lra & {\rm Ext}^1_H
     (\St_G ,\pi ) & \lra & \ds \bigoplus_{i=1,...,r} {\rm Ext}^1_H\,
     ({\rm ind}_{K_i}^H \, 1 , \pi )  \\
      & \cdots & \cdots & \cdots & \cdots & \cdots & \cdots \\
      & \lra  & {\rm Ext}^k_H\, (\St_H ,\pi ) & \lra &
     {\rm Ext}^k_H (\St_G ,\pi ) & \lra & \ds \bigoplus_{i=1,...,r} {\rm Ext}_H^k
     ({\rm ind}_{K_i}^H \, 1 , \pi ) \\
     & \cdots & \cdots & \cdots & \cdots & \cdots & \cdots \\
   \end{array}
   $$

   Let $i\in \{ 1,...,r\}$. Passing to contragredient representations, we obtain ${\rm Hom}_H\,
   ({\rm ind}_{K_{i}}^H 1, \pi ) \simeq {\rm Hom}_H \, ({\tilde \pi}, {\rm Ind}_{K_i}^H 1)$, the isomorphism
   depending functorially on $\pi$, and where $\rm Ind$ denote the functor of smooth induction. By
   Frobenius reciprocity for ${\rm Ind}$, we have ${\rm Hom}_H \, ({\tilde \pi}, {\rm Ind}_{K_i}^H 1)$ $\simeq$
     ${\rm Hom}_{K_i}({\tilde \pi}, 1)$, the isomorphism depending functorially on $\pi$. Passing to derived
     functor, we have ${\rm Ext}_H^k ({\rm ind}_{K_i}^H 1 , \pi )\simeq {\rm Ext}_{K_i}^k ({\tilde \pi}, 1)=0$, for
     all $k\geqslant 1$, since the group $K_i$ is compact. 

   \smallskip

   Hence our long exact sequence writes:
   
    $$
   \begin{array}{ccccccc}
     0 & \lra & {\rm Hom}_H\, (\St_H ,\pi ) & \lra & {\rm Hom}_H (\St_G ,\pi ) &
     \lra & \ds \bigoplus_{i=1,...,r} {\rm Hom}_{K_i} ({\tilde \pi},1)\\
     &  \lra & {\rm Ext}^1_H\, (\St_H ,\pi ) & \lra & {\rm Ext}^1_H
     (\St_G ,\pi ) & \lra & 0 \\
      & \cdots & \cdots & \cdots & \cdots & \cdots & \cdots \\
      & \lra  & {\rm Ext}^k_H\, (\St_H ,\pi ) & \lra &
     {\rm Ext}^k_H (\St_G ,\pi ) & \lra & 0 \\
     & \cdots & \cdots & \cdots & \cdots & \cdots & \cdots \\
   \end{array}
   $$

   Since $H$ is semisimple of $F$-rank $1$, the cohomological dimension of
   the category of finite length
   smooth complex representations of $H$ is $1$ (cf.  \cite{SS} {\S}III.3 or \cite{Be}{\S}4,  Theorem 29).
   So for $k\geqslant 2$, the  exact sequence
   $$
   0\lra   {\rm Ext}^k_H\, (\St_H ,\pi )  \lra 
   {\rm Ext}^k_H (\St_G ,\pi )  \lra  0
   $$
   \noi gives
   \begin{equation}
     {\rm Ext}^k_H \, (\St_G ,\pi )=0, \ \ k\geqslant 2
   \end{equation}

   We shall need the following result.

   \begin{proposition} Let $\sigma$ be an irreducible smooth complex
     representation of $H$. Then
     
     $${\rm Ext}^1_H \, (\St_H ,\sigma ) \simeq
     \left\{\begin{array}{cl}
     \CC & {\rm if} \ \ \sigma\simeq 1_H\\
     0 & {\rm otherwise}\end{array}\right.$$
     
   \end{proposition}

\noi {\it Proof}. This is a particular case of \cite{NP}, Theorem 2. 
\medskip

   Let us make three cases.

   \smallskip

   \noi {\it Case 1. Assume that $\pi$ is neither isomomorphic to the trivial
     representation of $H$ nor to $\St_H$}. Then by the previous proposition,
   the long exact sequence simplifies
   to give:
   $$
   0 \lra {\rm Hom}_H (\St_G ,\pi )\lra \ds \bigoplus_{i=1,...,r} {\rm Hom}_{K_i}({\tilde \pi}, 1) \lra 0
   \lra {\rm Ext}_H^1 (\St_G ,\pi )\lra 0\ .
   $$
   \noi Hence we obtain:
   
   \begin{equation}
     {\rm dim}_\CC \, {\rm Hom}_H \, (\St_G , \pi )=\sum_{i=1,...,r}
     {\rm dim}_\CC \, {\rm Hom}_{K_i}({\tilde \pi}, 1)
   \end{equation}
   \begin{equation}
     {\rm Ext}^1_H \, (\St_G ,\pi )=0
   \end{equation}

   \noi {\it Case 2. Assume that $\pi =\St_H$.} Then by Proposition 3.3, the
   long exact sequence writes:
   $$
   0\lra \CC \lra {\rm Hom}_H\, (\St_G ,\St_H )\lra \bigoplus_{i=1,...,r}
   {\rm Hom}_{K_i}(\St_H ,1) \lra 0\lra Ext_H^1 (\St_G ,\St_H )\lra 0
   $$
   
   It follows that
    
   \begin{equation}
     {\rm dim}_\CC \, {\rm Hom}_H \, (\St_G , \St_H )=\sum_{i=1,...,r}
      {\rm dim}_\CC \,  {\rm Hom}_{K_i}(\St_H ,1)  +1
   \end{equation}
   \begin{equation}
     {\rm Ext}^1_H\,  (\St_G ,\St_H )=0
   \end{equation}

   \noi {\it Case 3. Assume that $\pi =1_H$.} By Proposition 3.3, the long exact
   sequence writes
   $$
   0 \lra {\rm Hom}_H \, (\St_G ,1)\lra \bigoplus_{i=1,...,r} \CC
   \lra \CC \lra {\rm Ext}^1_H\, (\St_G ,1)\lra 0
   $$
   \noi This sequence does not allow us to compute ${\rm Hom}_H\, (\St_G ,1)$
   and ${\rm Ext}_H^1 (\St_G ,1)$ separatly. However we obtain:

   \begin{equation}
     {\rm EP}_H \, (\St_G ,1 )= r-1
   \end{equation}
   \begin{equation}
     {\rm dim}_\CC \, {\rm Hom}_H \, (\St_G ,1)\leqslant r , \ \ {\rm dim}_\CC
     \, {\rm Ext}_H^1\, (\St_G ,1)\leqslant 1
   \end{equation}

   Let us notice that  when $F$ is a finite field, we have
   ${\rm Ext}_H^1\, (\St_G ,1)=0$, so that
   \begin{equation}
     {\rm dim}_\CC\, {\rm Hom}_H (\St_G ,1)=r-1
   \end{equation}

   We may summarize our results as follows.

   \begin{theorem} \label{main} Let $\pi$ be an irreducible smooth
     representation of
     $H$.

     \noi a) If $\pi\not\simeq 1_H$, then ${\rm Ext}_H^k (\St_G ,\pi )=0$,
     for all $k\geqslant 1$. 

     \noi b) If $\pi\not\simeq 1_H$ then
     $$
          {\rm dim}_\CC \, (\St_G ,\pi )=
          \left\{\begin{array}{cc}
              {\ds  \sum_{i=1,...,r}} {\rm dim}_\CC \, {\rm Hom}_{K_i}({\tilde \pi}, 1) &
              {\rm if}\ \ \pi\not\simeq \St_H\\ 
              {\ds  \sum_{i=1,...,r}} {\rm dim}_\CC \,   {\rm Hom}_{K_i}({\tilde \pi}, 1)+1 & {\rm if}\ \ 
              \pi\simeq\St_H
          \end{array}\right.
          $$

          \noi c) We have ${\rm EP}_H \, (\St_G ,1) =r-1$.
     \end{theorem}

   \section{The case of ${\rm PGL}(2)$}

   In this section fix a quadratic extension $E/F$ and denote by $\theta$ the
   non-trivial element of  ${\rm Gal}(E/F)$. We may choose
   $\tau\in E\backslash F$ satisfying $\tau^2 =:\tau_0\in F$.
   We set $G={\rm PGL}(2,E)$ and
   still denote by $\theta$ the automorphism of $G$ obtained by applying
   $\theta$ to entries of matrices. Finally we set $H=G^\theta ={\rm PGL}(2,F)$.
   Then $H$ and $G$ are the sets of $F$-rational points of the reductive
   $F$-groups  $\HH ={\rm PGL}(2)$ and $\GG ={\rm Res}_{E/F}\, {\rm PGL}(2)$
   respectively. 
    
   \smallskip

   The set ${\mathcal P}_F (\GG )$ (resp. ${\mathcal P}_F (\HH )$) identifies
   with the projective line $P^1 (E)$ (resp. $P^1 (F)$). It is an easy exercise
   to prove that $H$ acts transitively on $P^1 (E) \backslash P^1 (F) =
   E\backslash F$, in the standard models of the projective lines $P^1 (E)
   =E\cup \{ \infty\}$, $P^1 (F)=F\cup \{ \infty\}$. It follows that
   $r=1$. With the notation of {\S}2 choose $P_1 ={\rm Stab}_G (\tau )$.
   Then $K_1 = {\rm Stab}_H (\tau )$ is the image in $H$ of the subgroup
   ${\tilde K}_1$ of ${\rm GL}(2,F)$ given by
   $$
   {\tilde K}_1 =\left\{ \matrice{a}{b\tau_0}{b}{a}\ ; \ a,b\in F, \ a^2 -
   b^2 \tau_0\not= 0\right\}
   $$

   In other words $K_1 \simeq E^\times /F^\times$. Irreducible representations of 
   $H$ distinguished by $E^\times /F^\times$ are known. More precisely, by
   \cite{W},  Lemme 9, page 219, if $\pi$ is a generic representation of $H$  
   then ${\rm dim}\, {\rm Hom}_{E^\times /F^\times} \, (1,\pi )\leqslant 1$, with
   equality when $\pi$ is not square integrable.

   It is an easy exercise that
   $E^\times /F^\times$ acts transitively on ${\mathcal P}_F (\HH )$ and that
   ${\mathcal P}_F (\HH )\simeq E^\times /F^\times$ as $E^\times /F^\times$-sets.
   Hence $\St_H$ is not $E^\times /F^\times$-distinguished. More
   generally, for a quadratique character $\chi$ of $F^\times$,
   $\chi\otimes\St_H$ is $E^\times /F^\times$-distinguished if and only if
   $\chi\circ N_{E/F}$ is non-trivial on $E^\times$, that is $\chi$ is
   non-trivial on $N_{E/F}(E^\times )$. In other words, $\chi\otimes \St_H$ is
   $E^\times /F^\times$-distinguished if, and only if,
   $\chi\not\in \{ 1, \eta_{E/F}\}$, where $\eta_{E/F}$ is the quadratic character
   of $E^\times$ attached to the extension $E/F$. 

   If  $\sigma$ is an irreducible  supercuspidal of $H$, the question of knowing whether $\pi$
   is $E^\times /F^\times$-distinguished is solved by a theorem of Tunnell's
   (see \cite{PrTunnell}, Theorem 1.1): the condition is given in terms of
   a value of the $\epsilon$-factor of the base change ${\rm BC}_{E/F}(\sigma )$ of
   $\sigma$ to $G$. 

   If $\pi =\chi\circ {\rm det}$, for some quadratic character
   $\chi$ of $F^\times$, then $\pi$ is  $E^\times /F^\times$-distinguished if,
   and only if, $\chi$ is trivial on $N_{E/F}(E^\times )$, that is
   $\chi \in \{ 1,\eta_{E/F}\}$.

   Finally, it is known (see e.g. \cite{Ha} Theorem 7.1) that $\St_G$ is not
   $H$-distinguished. Hence
   from ${\rm EP}\, (\St_G ,1)=r-1 =0$ and ${\rm Hom}_H \, (\St_G ,1)=0$, we
   obtain ${\rm Ext}^1_H \, (\St_G ,1)= 0$. In fact (loc. cit.), if $\chi$ is a quadratic
   character of $F^\times$, then ${\rm Hom}_H \, (\St_G ,\chi\circ {\rm det})$ is
   trivial if $\chi\not=\eta_{E/F}$ and $1$-dimensional if $\chi =\eta_{E/F}$. 

   \smallskip

   In the case of ${\rm PGL}(2)$, Theorem \ref{main} takes the following simplified
   form.

   \begin{theorem} a) For all irreducible smooth representations $\pi$ of
     ${\rm PGL}(2,F)$, we have
     $${\rm Ext}^1_{{\rm PGL}(2,F)}\, (\St_{{\rm PGL}(2,E)}, \pi )=0\ .
     $$
     \noi b) Let $\chi$ be a quadratic character of $F^\times$. We have
     $$
      {\rm dim}_\CC \, {\rm Hom}_{{\rm PGL}(2,F)}\, (\St_{{\rm PGL}(2,E)}, \chi\otimes
          \St_{{\rm PGL}(2,F)})=\left\{ \begin{array}{cc}
          0 & {\rm if}\ \ \chi=\eta_{E/F}\\
          1 & {\rm otherwise}\end{array}\right.
      $$

   \noi c) Let $\chi$ be a quadratic character of $F^\times$. We have
   
       $$
   {\rm dim}_\CC \, {\rm Hom}_{{\rm PGL}(2,F)}\, (\St_{{\rm PGL}(2,E)},
   \chi\circ {\rm det})=\left\{ \begin{array}{cc}
          1 & {\rm if}\ \ \chi=\eta_{E/F}\\
          0 & {\rm otherwise}\end{array}\right.
       $$
   \noi d) If $\pi$ is a principal series representation of ${\rm PGL}(2,F)$
   then
   $$  {\rm dim}_\CC \, {\rm Hom}_{{\rm PGL}(2,F)}\, (\St_{{\rm PGL}(2,E)},\pi )=1
   $$

   \noi e) If $\pi$ is an irreducible supercuspidal representation of
   ${\rm PGL}(2,F)$, then
   $$ {\rm dim}_\CC \, {\rm Hom}_{{\rm PGL}(2,F)}\, (\St_{{\rm PGL}(2,E)},\pi )=
   {\rm dim}_\CC \, {\rm Hom}_{E^\times /F^\times} ({\tilde \pi}, 1) \in \{ 0,1\}
   $$
\end{theorem}

\section{Some comments and questions}

\noi {\bf On the number $r$}. Recall that $r$ is the number of orbits of $H$
in the set of minimal $\theta$-split parabolic $F$-subgroups of $\GG$.

\begin{proposition} When $F$ is finite, we have $r\leqslant 2$.
\end{proposition}

\noi {\it Proof}. Since any minimal $\theta$-split parabolic $F$-subgroup of
$\GG$  contains a maximal $(\theta ,F)$-split torus of $\GG$, and since
any maximal $(\theta ,F)$-split torus is contained in exactly two
minimal $\theta$-split parabolic $F$-subgroup ($\GG$ has $F$-rank $1$),
we are reduced to proving that
any two maximal $(\theta ,F)$-split tori  are conjugate in $H=\HH^\circ (F)$.

So let $A_1$, $A_2$ be maximal $(\theta ,F)$-split tori of $\GG$. By \cite{HW}
Proposition 10.3, there exists $g\in (\HH^\circ Z_G (A_1 ))(F)$ such that
$A_2 =gA_1g^{-1}$.  Since the Levi subgroup ${\mathbb L}=Z_G (A_1 )$ is
connected, so is $\HH^\circ \cap {\mathbb L}$. Hence by Lang's Theorem
the first Galois cohomology set $H^1 (F, \HH^\circ\cap {\mathbb L})$ is
trivial. Now by a classical cocycle argument, we have that
$$
(\HH^\circ Z_G (A_1 ))(F) = \HH^\circ (F) \, Z_G (A_1 )(F) = H Z_G (A_1 )(F)\
$$

\noi and our result follows.

\begin{proposition}\label{Galois} Assume that $F$ is local and that $\GG /\HH$ is a Galois
  symmetric space. Then $r\leqslant 2$.
\end{proposition}

\noi {\it Proof}. By assumption $\HH$ is a connected semisimple $F$-group of $F$-rank $1$,
$\GG = {\rm Res}_{E/F}\, \HH$, for some quadratic extension $E/F$, and $\GG$ is assumed to
have $E$-rank $1$.  As a consequence of  a conjecture of D. Prasad's proved by R.
Beuzart-Plessis (\cite{BP}),  have:
$$
{\rm dim}_\CC\, {\rm Hom}_H \, (\St_G ,1 )\leqslant 1
$$
\noi On the other hand, by Theorem 3.4, we have
$$
     {\rm EP}_H \, (\St_G ,1)=r-1
$$
 \noi It follows that $r-1\leqslant 1$, that is $r\leqslant 2$.  
 \smallskip

 In all cases where I could manage to compute $r$, I always found out that
 $r\leqslant 2$. This is indeed the case in examples  $1$, $2$ and $4$ of {\S}1.
 It is natural to ask whether this inequality holds true in general. 
 \bigskip

 \noi {\bf On vanishing of ${\rm Ext}$ groups}. With the notation of Theorem 3.4, we have
 that ${\rm Ext}_H^k (\St_G ,\pi )$ vanishes for all irreducible smooth representations of $H$
 not isomorphic to the trivial representation, and all $k\geqslant 1$. The question of whether
 ${\rm Ext}_H^1 (\St_G ,1_H )$ vanishes or not in general remains open.

 \bigskip

 \noi {\bf On the quotients $H/K_i$}. A striking consequence of Theorem \ref{main} is that
 the branching  problem for the Steinberg representation of $G$ is related to distinction problems for  
 related to ``smaller  quotients''
 $\HH^\circ /{\mathbb K}_i$, $i=1,...,r$. It is therefore natural to ask whether the
 homogeneous spaces $\HH^\circ /{\mathbb K}_i$ are reductive symmetric spaces or not.
 I all cases where I can push computations to the end I never found any counter-example.

 Let us however quote the following partial result.

 \begin{proposition} \label{symspace} Assume that the relative root system of $G$ is reduced. Then
   the quotients $\HH^\circ /{\mathbb K}_i$, $i=1,...,r$, are reductive symmetric spaces.
 \end{proposition}

 \noi {\it Proof}.  Write $\HH^\circ /\KK$ for one of the quotients $\HH^\circ /\KK_i$, $i=1,...,r$.
 There exists a maximal $(\theta ,F)$-split torus $T$ in $\GG$ such that $\KK = \HH^\circ\cap\LL$, where
 $\LL$ is the Levi subgroup $Z_{\GG}(T)$. By definition, we must prove that there exists a rational involution
 of $\HH^\circ$ with connected fixator $\KK^\circ$.
 \smallskip

 \noi {\it Case $1$. Assume that the center $Z(\GG )$ of $\GG$ is trivial}. Since $T$ is
 $F$-isomorphic to the multiplicative group $\GG_m$, there exists $t_0 \in T(F)$, such that
 $t_0 \not= 1$ and $t_0^2 =1$. In particular $t_0\not\in Z(G)$. Set $\tau ={\rm Int}(t_0 )$;
 this is a non trivial rational involution of $\GG$.

 Let $\{ \pm \alpha \} = \Phi (T,\GG )$ be the relative root system of $T$ in $\GG$. Since
 $t_0$ is not central, we have $\alpha (t_0 )\not= 1$. It follows that $t_0$ is a regular
 element of $T$ and that $Z_\GG (t_0 )^\circ =\LL$\footnote{Note that if $\Phi (T,\GG )$
   is not reduced, then $Z_\GG (t_0 )^\circ$ is bigger that $\LL$.}. In other words $(\GG^\tau )^\circ = \LL$. 

 Since $\theta$ stabilizes $T$, we have $\theta (\LL )=\LL$. Consider the two rational involutions
 of $\GG$ given by $\tau$ and  $\theta\tau\theta^{-1}$. The connected fixator of the first one is
 $\LL$ and the connected fixator of the second one is $\theta (\LL )=\LL$.  By \cite{HW} Proposition
 1.2, a rational involution of  $\GG$ is entirely determined by its connected fixator. It follows
 that $\theta\tau\theta^{-1}=\tau$, that is: $\tau$ and $\theta$ commute.

 In particular $\tau$ stabilizes $\HH^\circ$ and ${\tilde \tau}=\tau_{\mid \HH^\circ}$ is a rational
 involution of
 $\HH^\circ$. We have $((\HH^\circ )^{\tilde \tau})^\circ \subset \HH^\circ \cap (\GG^\tau )^\circ$, that
 is $((\HH^\circ )^{\tilde \tau})^\circ \subset \HH^\circ \cap \LL =\KK$. On the other hand the
 containment $\KK^\circ  \subset ((\HH^\circ )^{\tilde \tau})^\circ$ is straightforward. So we
 finally obtain
 $\KK^\circ  = ((\HH^\circ )^{\tilde \tau})^\circ$, as required. 
 \medskip

 \noi {\it Case $2$. The center of $\GG$ is not necessarily trivial}. Abreviate $Z=Z(\GG )$; this is
 a finite group.  Set $\bGG =\GG /Z$ and $\bLL =\LL /Z\subset \bGG$. Observe that $\bGG$ has
 trivial center.  Denote by $\bT$ the image of $T$ in $\bGG$. This is a maximal $F$-split torus of
 $\bGG$.

 Let us first prove that the centralizer of $\bT$ in $\bGG$ is $\bLL$. Let $g\in \GG$. Then $gZ$
 lies in the centralizer of $\bT$ if and only if $gtg^{-1} \in tZ$ for all $t\in T$. If
 $n=\vert Z\vert$ denotes the order of $Z$, we obtain $gt^n g^{-1} =t^n$ for all $t\in T$.
 Since the map $T\lra T$, $t\mapsto t^n$,  is surjective, we have that $g\in Z_\GG (T)=\LL$, as
 required.

 Applying the construction of Case $1$, there is an element ${\bar t}_0\in \bT$ such that ${\bar t}_0
 \not= 1$, ${\bar t}_0^2 =1$ and $\bLL =(Z_{\bGG}({\bar t}_0))^\circ$. Write $\bt_0 =t_0 Z$, for some
 $t_0 \in T$. Note that it is not always possible to choose $t_0$ in $T(F)$. By construction
 $t_0\not\in Z$ and $t_0^2\in Z$. So $\tau := {\rm Int}(t_0 )$ is a non trivial  involution of
 $\GG$. Let us prove that $Z_\GG (t_0 )^\circ =\LL$. Because $Z\subset \LL$ and $\LL$ is
 connected, one containment is obvious. Let $p$~: $\GG\lra \bGG$ be the projection. One easily sees
 that $p(Z_\GG (t_0 ))=Z_{\bGG} (\bt_0 )$. Hence
 $p(Z_\GG (t_0 )^\circ )\subset Z_\GG (\bt_0 )^\circ =\LL /Z$, and we are done.

 Since the connected fixator of $\tau$ is defined over $F$, it follows from \cite{HW} Proposition
 1.6 that $\tau$ itself is defined over $F$.

 As in case $1$, we prove that $\theta$ and $\tau$ commute and that ${\tilde \tau}=\tau_{\mid \HH^\circ}$
 is a rational involution of $\HH^\circ$ with fixator $\KK^\circ$. 
 \medskip

 \noi {\bf An example where multiplicity $1$ fails}. The following example answers a question
 of Waldspurger's.
 \smallskip

 For $G/H$ take the quotient ${\rm SL}(2,E)/{\rm SL}(2,F)$, where $E$ is a quadratic extension of our field
 $F$ assumed to have odd residue characteristic. In that case
 we have ${\rm dim}\, {\rm Hom}_H \, (\St_G ,\St_H )\geqslant 3$\footnote{In fact we have equality:
   ${\rm dim}\, {\rm Hom}_H \, (\St_G ,\St_H ) =  3$. To prove this one has to use the model of the
   algebraic dual of the Steinberg representation given in terms of  harmonic cochains as in \cite{Br}.}
 \smallskip
 
 Identify $H$ with ${\rm SL}_F (E)$. Let $N=N_{E/F}$~: $E^\times\lra F^\times$ denote the norm map and set
 $E^1 ={\rm ker}\, N$
 that we identify with a subgroup of $H$. Here $r=2$ and the groups $K_i$, $i=1,2$, are conjugate under
 ${\rm GL}(2,F)$ to  $E^1$. By Theorem \ref{main} and since $\St_H$ is selfdual, we have :
 $$
 {\rm dim}\, {\rm Hom}_H \, (\St_G ,\St_H ) = 2\times {\rm dim}\, {\rm Hom}_{E^1}\, (\St_H ,1 ) +1
 $$

 \noi We are going to prove that ${\rm dim}\, {\rm Hom}_{E^1}\, (\St_H ,1 ) \geqslant 1$; our result will follow.
 Observe that ${\rm Hom}_{E^1}\, (\St_H ,1 ) ={\rm Hom}_\CC (\St_H ,1 )^{E^1}$, the space of $E^1$-invariant
 linear forms on $\St_H$. It contains $({\tilde \St}_H )^{E^1}$, where ${\tilde \St}_H$ denotes the space of smooth
   linear forms on $\St_H$, that is the contragredient of $\St_H$, which turns out to be $H$-isomorphic to
   $\St_H$, $\St_H$ being self-dual. To sum up we have
   $$
   {\rm Hom}_{E^1}\, (\St_H ,1 ) \supset \St_H^{E^1}
   $$
   \noi We are going to prove that ${\rm dim}\, \St_H^{E^1}=1$ and our claim will follow.

 Recall that $\St_H =V/V_0$, where $V$ is the space of locally constant functions on $P^1 (F)$ and
 $V_0$ the subspace of constant functions.  The group $E^1$ has two orbits in $P^1 (F)$. Indeed one may identify
 $P^1 (F)$ with $E^\times /F^\times$ so that the quotient space $E^1 \backslash P^1 (F)$ is in bijection
 with $E^\times /(E^1 F^\times )$. The norm maps induces an exact sequence
 $$
 1\lra E^1 F^\times \lra E^\times \lra N(E^\times )/(F^\times )^2 \lra 1
 $$
 \noi Hence
 $$
 \vert E^\times /E^1 F^\times\vert =\vert N(E^\times )/(F^\times )^2 \vert=
 \frac{\vert F^\times /(F^\times )^2\vert}{\vert F^\times /N(E^\times )\vert} = 4/2 =2
 $$
 \noi where equality $\vert F^\times /N(E^\times )\vert =2$ follows from class field theory for
 quadratic extensions.

 Let $\Omega_i$, $i=1,2$ be the orbits of $E^1$ in $P^1 (F)$. Let $e\in E^\times \backslash E^1 F^\times$. Then the
 map $P^1 (F)\lra P^1 (F)$,  mapping a line $d$ to $ed$,  is a homeomorphism taking $\Omega^1$ to $\Omega_2$.
 Since at least one of the $\Omega_i$'s is open, they are both open. From this it follows that the fixed
 space $V^{E^1}$ is two dimensional, whence that
 $\St_H^{E^1} =V^{E^1}/V_0$ is $1$ dimensional,  as we claimed.

paul.broussous{@}math.univ-poitiers.fr
\medskip

D\'epartement de Math\'ematiques

UMR 7348 du CNRS
\smallskip

T\'el\'eport 2 - BP 30179

Boulevard Marie et Pierre Curie

86962 Futuroscope Chasseneuil Cedex

 France

\end{document}